\documentclass[11pt, a4paper]{article}
\usepackage{amsfonts}
\usepackage{float}
\usepackage{amssymb}
\usepackage{epsfig}
\usepackage{caption}
\usepackage[table]{xcolor}
\usepackage{url}
\usepackage{amssymb}
\usepackage{pifont,booktabs}
\usepackage{lineno}
\newcommand{\cmark}{\ding{51}}
\newcommand{\xmark}{\ding{55}}

\usepackage{pgfplots}
\usepackage{subfigure}
\usepackage{array}
\newcolumntype{C}[1]{>{\centering\arraybackslash}p{#1}}

\newcommand{\vv}[1]{\underline{#1}}

\newcommand{\mod}[1]{\left\vert{#1}\right\vert}
\newcommand{\norm}[1]{\left\Vert{#1}\right\Vert}
\newcommand{\en}[1]{(\ref{#1})}
\newcommand{\rr}{\mathbb{R}}

\newcommand{\corr}[1]{\color{black}{#1}}

\newcommand{\apap}[1]{``{#1}}
\newcommand{\apch}{"}

\newtheorem{teor}{Theorem}

\title{Results, challenges and new steps on RBF approximation and computation\thanks{Work partially supported by INdAM-GNCS research funds.}}

\author{Stefano De Marchi\\ Department of Medicine\\ University of Padova, Italy
\and Nadaniela Egidi, Josephin Giacomini, Pierluigi Maponi\\
School of Science and Technology\\ University of Camerino, Italy}

\date{}
\begin{document}

\maketitle

\begin{abstract}
We present an up-to-date overview of approximation methods based on Radial Basis Function (RBF) techniques, which have recently attracted attention across various computational tools and application fields. This review study presents relevant results on RBF techniques, highlighting the associated computational challenges and stability issues that must be addressed when high-performance or parallel computation is required. 
\end{abstract}

\section{Introduction}
\label{sec1}
The origins of RBF (Radial Basis Function) techniques date back to the second half of the last century and stem from a geodetic application by Hardy \cite{hardy}, in which irregular surfaces in cartography were obtained from sparse topographical data. In this seminal work, the representation of complex surfaces was proposed as a sum of quadric forms (such as cones or hyperboloids) centered at different data points. The simple idea of approximating data at scattered centers attracted the attention of several investigators, leading to the mathematical foundation of RBF or, more generally, to kernel-based approximation. The first relevant results about convergence, stability, and existence of solutions were provided in \cite{Meinguet}. In \cite{Franke}, a critical benchmark test of various scattered data interpolation methods was presented, showing that the multiquadric method yielded the most accurate results. Another fundamental theoretical result is the proof for the solvability of the RBF interpolation matrix \cite{Micchelli}. Furthermore, the theory of  RBF interpolation can rely on a physical interpretation since the family of the so-called {\em Thin Plate Spline} can be obtained as the unique minimizer of a specific semi-norm (representing \apap bending energy") \cite{Duchon}. 

These fundamental results provide a suitable computational environment for addressing various problems. RBFs are primarily valued for their mesh-free nature, which allows them to handle complex geometries and scattered data without a structured grid. For example, the Kansa method \cite{Kansa} allows one to compute solutions to Partial Differential Equations (PDEs) without a rigid geometric grid, unlike finite element approaches.
RBFs are now {{reference tools}} for mesh morphing, in which an existing computational grid (e.g., around an aircraft wing) is smoothly deformed using a few control points while preserving the quality of the cells in the internal grid (cf. \cite{Biancolini, Boer}). The resulting deformation is smooth, mesh-independent, and preserves element quality, making RBFs particularly attractive for shape optimization, fluid–structure interaction, and moving-boundary problems.

RBFs are also used for elastic registration \cite{Rohr}, in which two different images (e.g., an MRI - Magnetic Resonance Imaging - and a CT - Computed Tomography - scan of the same patient) are deformed to perfectly align the anatomical landmarks.
This {\em meshless} approach has indeed become quite standard for complex engineering and medical simulations.

RBFs also play a crucial role in scientific Machine Learning, both in Artificial Neural Networks (ANNs) and as kernels in Support Vector Machines (SVMs). 
The close connection between training ANNs and data-fitting processes is discussed in \cite{Broomhead} and addressed using an RBF approach. 
Several contributions further strengthened this initial work on RBF networks. For instance, in \cite{Moody} the authors demonstrated the efficiency of RBFs in learning complex mappings through locally tuned hidden units.
The connection between RBFs and ANNs was generalized in \cite{Poggio}, where RBF networks are shown to be universal approximators. 
The kernel trick was initially proposed in \cite{Boser} to integrate kernel functions into the SVM framework. It demonstrated that SVMs could use RBF kernels to automatically select \apap supporting patterns" (support vectors) as the centers of the radial functions, allowing the model to adapt its complexity to the data.

Although RBFs are mathematically elegant approximation tools, they face several critical technical challenges that can compromise their effectiveness. The {\em Shape Parameter dilemma} affects most usual kernels, like Multiquadrics or Gaussians, and relies on the practical choice of the parameter that controls the basis function's width. Its choice is therefore crucial for obtaining effective implementations of RBF tools \cite{Fasshauer}.
The {\em numerical instability} is strictly related to the shape parameter and the ill-conditioning of the interpolation matrix as the number of data points increases \cite{Fornberg}. 
The {\em computational complexity and scalability} are other critical aspects of RBF techniques. In fact, 
differently from what happens in Finite Element Methods (FEMs), classic RBF matrices are dense when global RBFs are used, so the solution of an RBF interpolation problem scales at $O(N^3)$ for the computational load and at $O(N^2)$ for the storage load, where $N$ is the characteristic dimension of the problem.
In addition, the {\em selection of centers} of the RBFs is a non-trivial optimization problem \cite{DeMarchi, DeMarchi2, Wenzel} that affects the accuracy of the overall approximation process. The problem is also {\corr{related}} to the choice of a stable basis, as discussed in \cite{
Muller2}.

All these critical aspects must be properly addressed when the problem to be solved is large and/or difficult due to the problem's intrinsic instabilities. Thus, they are the key to making RBF a relevant tool in High Performance Scientific Computing (HPSC) \cite{Eijkhout}.  In the following sections, we provide an overview of research efforts addressing critical aspects of RBF by selecting a few works. The selected material is also organized into two directions: those that mainly address the computational issues of RBF and those that address its stability issues. Of course, these aspects are never separable in any computational procedure, but this organization allows us to focus on the main contributions of the presented works. 

The paper is organized as follows. In Section \ref{sec2}, we present the main achievements in the stabilization of RBF techniques. In Section \ref{sec3}, we discuss the results for efficient computer implementations of RBF techniques. In Section \ref{sec4}, we outline some conclusions and final remarks.

\section{Conditioning and stabilization strategies}
\label{sec2}
We start by recalling the RBF interpolation problem and introducing in Table \ref{tab:notation} the notation that will be used in the following.

\begin{table*}[ht]
\caption{Notations and symbols.}
\label{tab:notation}
\centering
{\footnotesize 
\begin{tabular}{p{1.3cm}p{4.0cm}||p{1.3cm}p{4.0cm}}
\toprule 
\textbf{Notation} &
\textbf{Description} &\textbf{Notation} &
\textbf{Description} \\
\midrule

$D$
& domain, subset of $\rr^d$ &$\phi$
&
Radial basis function
\\
$X\subset D$
&
Set of centers & $\phi_j$
&
RBF centered at $\vv x_j\in X$
\\
$P_f$
&
Interpolator of $f$ & $N_\phi(D)$
&
Native space
\\
$h_{X,D}$
&
Fill distance &$P_{\phi,X}$
&
Power function
\\
$q_{X}$
&
Separation distance &$\|\cdot \|$ & 2-norm
\\
$\rho$
&
Uniformity & $\|\cdot \|_{N_\phi}$ & Native space norm
\\
\bottomrule
\end{tabular}
}
\end{table*}

Let $D\subset\rr^d$ some subset of $\rr^d$. Take the set of points $(\vv x_i,f_i)$, $i=1,2,\ldots,N,$ being the interpolation data, where $X=\{\vv x_i\in D\subset\rr^d,\ i=1,2,\ldots,N\}$ is the set of the centers (often corresponding to the interpolation points) and $f_i\in\rr$, $i=1,2,\ldots,N,$ are the interpolation values. 

The interpolation problem is: given functions $\phi_i:D\to\rr$, $i=1,2,\ldots,N$, compute the coefficients $c_j$, $j=1,2,\ldots,N,$ such that 
\begin{equation}
\sum_{j=1}^N c_j\phi_j(\vv x_i)=f_i,\quad i=1,2,\ldots,N,
\label{RBFint}
\end{equation}
of course, it is required that this solution is unique.
\begin{table}[h!]
\centering
\begin{tabular}{l|c|c|c|c}
Name & smoothness &$\phi$ & support & SPD\\
\hline
Gaussian  & $C^\infty$ &$e^{-r^2}$ & global & \cmark\\
IM& $C^\infty$   &$(1+r^2)^{-1/2}$ &global & \cmark\\
Mat\'ern & $C^2$ &$(1+r)e^{-r}$& global & \cmark\\
Wendland &$C^2$&$\max\{1-r,0\}^4(4 r+1)$& compact & \cmark\\
TPS&  $C^\infty$   &$r^2\log(r)$ & global & \xmark
\end{tabular}
\caption{Most used kernel functions in RBF approximation; the acronyms IM and TPS stand for Inverse Multiquadric and Thin Plate Spline, respectively. We note that the last column reports whether they are shape parameter dependent (SPD).}
\label{t1}
\end{table}

\begin{figure}[h!]
\small
\centering
\includegraphics[width=.96\textwidth]{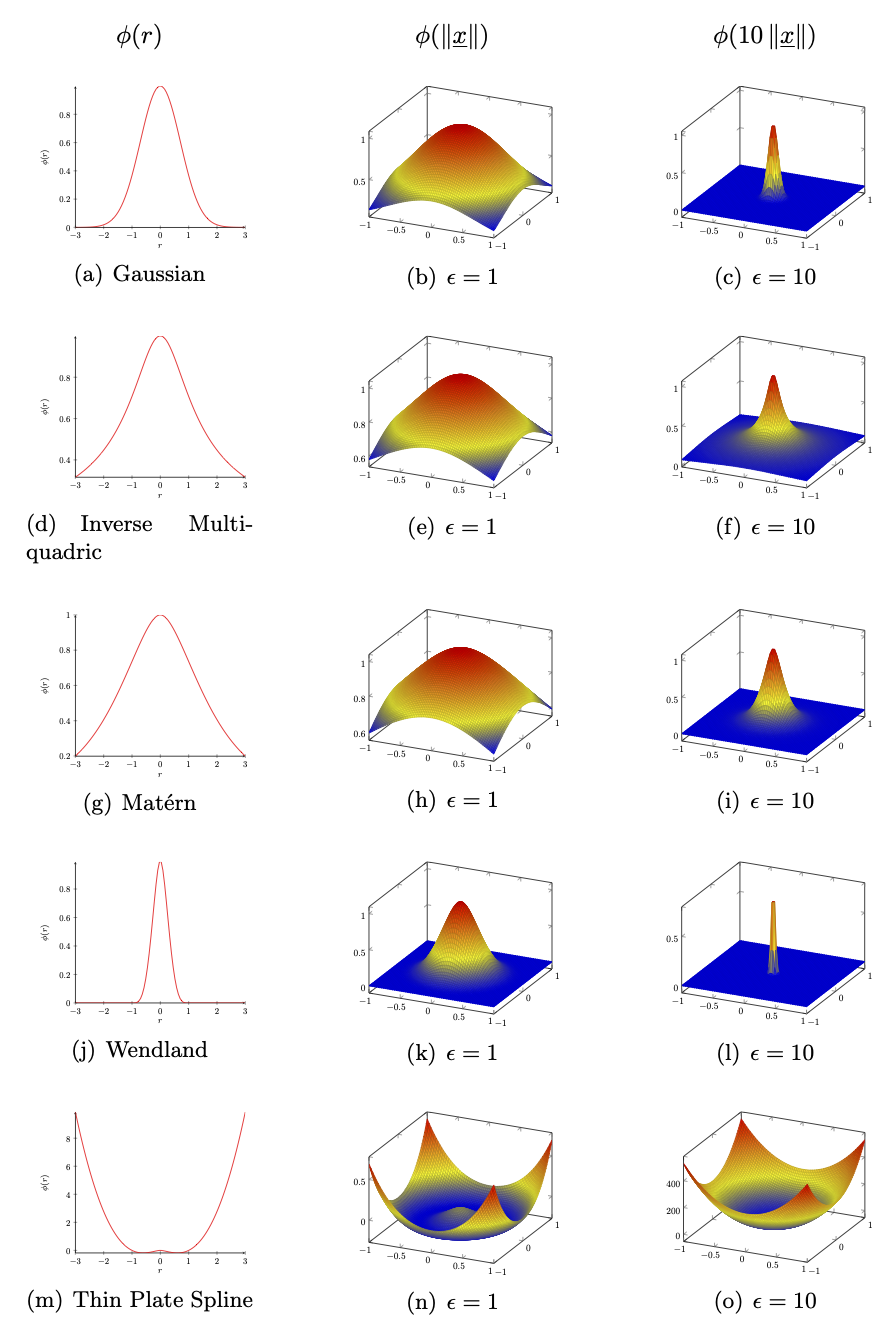}
\caption{Graphs of the kernels reported in Table \ref{t1}. For $\rr^2$ we used the shape parameters $\epsilon=1,\,10$.}
\label{f2}
\end{figure}

In the case of RBF, the functions $\phi_j$, $j=1,2,\ldots,N,$ are generated by a unique kernel function $\phi$, more precisely $\phi_j(\vv x)=\phi(\norm{\vv x-\vv x_j})$, $j=1,2,\ldots,N$, $\vv x\in D$, with $\phi:[0,\infty) \rightarrow \rr$.
Among the various RBFs available in the literature, we have reported in Table \ref{t1} the most frequently used in applications: the Gaussian, the Inverse Multiquadric (IM), the Mat\'ern, the Wendland, and the Thin Plate Splines (TPS). We note that the Mat\'ern and Wendland functions reported in the table are twice-differentiable members of their respective families. Also, TPS can be seen as a special case of the polyharmonic function family (cf. e.g. \cite{Fasshauer} for details). The shape of these functions, except for the TPS, can be easily controlled by the so-called {\it shape parameter} $\epsilon>0$ through the scaling of the kernel argument, that is, considering $\phi(\epsilon r)$. In Figure \ref{f2}, we have reported the kernel graphs (extended by parity) described in Table \ref{t1}, and the graphs of $\phi(\epsilon \norm{\vv x})$, $\vv x\in\rr^2$ for different choices of the shape parameter $\epsilon$. 

Assuming that $\phi$ is strictly positive definite (see e.g. \cite{Fasshauer}), the interpolation problem \en{RBFint} has a unique solution since the interpolation matrix $\Phi$, where $[\Phi]_{i,j}=\phi(\norm{\vv x_i-\vv x_j})\in\rr^{N\times N}$, 
is strictly positive definite. 
We note that this does not affect the generality of the discussion, since, by the Paley–Wiener Theorem, for every conditionally positive definite function of order $k$ (see e.g. \cite{Fasshauer}), there exists an associated strictly positive definite function (see e.g. \cite{DeMarchi2}). 
However, the linear system \en{RBFint} is usually ill-conditioned as $N$ grows. Several authors have investigated the spectral properties of the interpolation matrix $\Phi$; see, for example, \cite{Narcowich1, Narcowich2, Narcowich3, Ball1, Ball2, Braun}. In particular, \cite{Schaback} proved that the minimum eigenvalue $\lambda$ of $\Phi$ is bounded as follows:
\begin{equation}
cN^{-\beta/d}\le\lambda\le CN^{-\beta/d}
\label{eigbounds}
\end{equation}
where $c$, $C>0$ are constants depending on the interpolation points and the kernel function $\phi$, $\beta>0$ is a smoothness parameter for $\phi$. More precisely, the Fourier transform $\hat\phi$ satisfies the following relation:
\begin{equation}
0<\hat\phi(\omega)<A(1+\mod\omega)^{-\beta},\quad \omega\in \rr,
\label{smooth}
\end{equation}
for a suitable constant $A>0$. Furthermore, for smooth kernels, such as Gaussian and multiquadratic ones,
\en{eigbounds} holds with bounds that decay exponentially.

Obtained the $c_j$ as the solution of \en{RBFint}, we can consider the function
\begin{equation}
P_f(\vv x)=\sum_{j=1}^N c_j\phi_j(\vv x),\quad\vv x\in D,
\label{Pf}
\end{equation}
as an approximation for the function $f(\vv x)$, $\vv x\in D$ that generated the interpolation data $(\vv x_i,f_i)$, $i=1,2,\ldots,N$, i.e., $f(\vv x_i)=f_i$, $i=1,2,\ldots,N$. The accuracy of $P_f$ depends on the density of the point set $X$ in $D$, and on the regularity of the function $f$. The {\it fill distance} is the classical measure of the density of the interpolation points, it is given by
\begin{equation}
h_{X,D}:=\sup_{\vv x\in D}\min\left\{\norm{\vv x-\vv x_j},\  j=1,2,\ldots,N\right\},
\label{Fdistance}
\end{equation}
and represents the radius of the largest ball that can be placed among the interpolation points. 

Other quantities often used to analyze stability are the {\it separation distance}
\begin{equation}
    q_X:=\frac12\min\left\{\norm{\vv x_i-\vv x_j},\ 1\le i,j\le N, \ i\ne j\right\},
\end{equation}
and the {\it uniformity}
\begin{equation}
\rho:=\frac{q_X}{h_{X,D}}.
\end{equation}
In particular, if $\rho$ is near $1/\sqrt{d}$,
the interpolation data are nearly equally distributed in the Euclidean norm on $\rr^d$. Moreover, in contrast to polynomial interpolation, radial basis interpolants perform better when the interpolation nodes are nearly uniformly distributed \cite{DeMar}.

Classical error estimates are usually given for functions in $N_\phi(D)$, the {\em native space} of $\phi$ which is obtained from the infinite dimensional space $H_\phi(D)=span\{\phi(\norm{{\vv x - \cdot}}),\ \vv x\in D\}$ equipped with the scalar product:
\begin{equation}
\left(\sum_ip_i\phi(\norm{{\vv x_i-\cdot}}),\sum_jq_j\phi(\norm{{\vv x_j-\cdot}})\right)=\sum_{i,j}p_iq_j\phi(\norm{\vv x_i-\vv x_j}),
\label{Sprod}
\end{equation}
where finite or infinite summations can be considered. In particular, $N_\phi(D)$ is defined as the completion of $H_\phi(D)$ with respect to the norm $\norm{\cdot}_\phi$ induced by \en{Sprod}. We recall that the native space is the Reproducing Kernel Hilbert Space (RKHS) associated to $\phi$.

\begin{teor} \label{Perror}
Let $\phi: D\to\rr$ be a continuous function and a strictly positive definite kernel. Let $\vv x_i\in D\subset\rr^d$, $i=1,2,\ldots,N,$ be pairwise distinct points, and let $f\in N_\phi(D)$. Then for every $\vv x\in D$ we have:
\begin{equation}
\mod{f(\vv x)-P_f(\vv x)}\le P_{\phi,X}(\vv x)\norm{f}_{N_\phi},\quad\vv x\in D,
\label{error}
\end{equation}
where $P_{\phi,X}$ is the {\tt power function} which depends on the interpolation points $X$ and on the kernel $\phi$, and $\norm{f}_{N_\phi}$ is the native space norm of $f$. 
\end{teor}
See \cite[Chap. 14]{Fasshauer} for a detailed proof of this theorem.
\vskip 0.2in
Theorem \ref{Perror} assumes knowledge of the exact solution $c_j$, $j=1,2,\ldots,N$ of \en{RBFint}. However, only an approximated solution $\tilde c_j$, $j=1,2,\ldots,N$ of \en{RBFint} is usually known {{that yields}} the following error inequality 
\begin{equation}
\mod{f(\vv x)-\tilde P_f(\vv x)}\le \mod{f(\vv x)- P_f(\vv x)}+\mod{P_f(\vv x)-\tilde P_f(\vv x)},\quad\vv x\in D,
\label{error2}
\end{equation}
where $\tilde P_f$ is the function \en{Pf} {{calculated using the approximate}} solution $\tilde c_j$, $j=1,2,\ldots,N$, of \en{RBFint}. Due to the ill-conditioning of \en{RBFint}, the second addendum, in the right-hand side of \en{error2}, may add a non-negligible quantity to the interpolation error. 

Among the most popular methods for stabilizing the interpolant, many belong to the following categories.  
\begin{enumerate}
	\item \textbf{RBF-QR method}: it is rooted in a particular decomposition of the kernel, and it has been developed so far to treat the Gaussian kernel \cite{FoPi07,Fornberg11,FasshauerMcCourt12,Kormann19}.
	\item \textbf{Hilbert-Schmidt Singular Value Decomposition (HS-SVD)}: it has been developed to stably compute the RBF interpolants \cite{Cavoretto15bb,FasshauerMcCourt}. In principle, this technique can be applied to any kernel, provided that the HS eigenvalues and eigenvectors are known. However, these quantities are {{far from easy to compute and,}} in practice, they only work for the Gaussian function.
	\item \textbf{Weighted Singular Value Decomposition (WSVD) bases}: it is a more general approach that applies to any RBF, consisting of computing a weighted SVD {{that}} produces stable bases \cite{DeMarchi13}. 
\end{enumerate}
We recall that when we have flat or nearly flat RBFs, a severe ill-conditioning of the RBF interpolation matrix occurs, even if 
it generally {{guaranties}} the best accuracy for smooth solutions. The {\it RBF-QR method} tries to reduce the error due to the ill-conditioning when the shape parameter approaches zero, so the RBF becomes flat. In this case, the error term \en{error} is usually small, but the second addendum, in the right-hand side of \en{error2}, is usually large as a consequence of the ill-conditioning of \en{RBFint}. The core idea of RBF-QR is to change the basis of the approximation space to eliminate numerical instability; this is done by expanding the RBF kernel into a series of more stable functions, such as spherical harmonics \cite{Fornberg1} or Chebyshev polynomials \cite{Fornberg2}. A QR decomposition of the representation matrix arising from these expansions can gather the ill-conditioning in simple matrices that can be inverted almost analytically. A similar idea is proposed in \cite{Fasshauer2} for the Gaussian kernel, for which a Mercer expansion is considered together with the QR decomposition of the sampled eigenfunctions used in the expansion. 

{{Below we illustrate two other categories of stabilization techniques  of RBF interpolation: greedy techniques and Partition of Unity techniques.}}

Greedy techniques \cite{Schaback2} are iterative algorithms {{to select}} an optimal subset of points (centers)
from a large {{data set}} $X$ to construct an efficient and stable interpolant. 
Several variants have been proposed. 
The {\em p-greedy} \cite{
DeMarchi} is a data-independent method that iteratively selects the next point where the power function $P_{\phi, X}$ is maximized, ensuring in this way that the points are \apap well-spread" and minimize the worst-case error for any function in $N_\phi(D)$. 


\!\!\!\!\!\!\!\!\!\!\!\!\!\!\!\!\begin{table*}[ht]
\caption{Comparison of the main stabilization and complexity-reduction techniques for kernel-based interpolation.}
\label{tab:comparison_stabilization}
{\footnotesize
\!\!\!\!\!\!\!\!\!\!\!\!\!\!\!\!\begin{tabular}{p{1.3cm}p{2.0cm}p{2.0cm}p{1.3cm}p{1.3cm}p{1.4cm}p{2.2cm}}
\toprule
\textbf{Method} &
\textbf{Main idea} &
\textbf{Pros} &
\textbf{Cost} &
\textbf{Robust w.r.t. $\varepsilon$} &
\textbf{Scalability} &
\textbf{Cons} \\
\midrule

\rowcolor{blue!20} RBF-QR
&
Stable basis transformation
&
{\small High accuracy with Gaussian or MQ kernels in small-to-medium problems}
&
$O(N^3)$
&
Excellent
&
Moderate
&
{\small Kernel-specific implementation; difficult extension to high dimensions.} \\

\rowcolor{blue!20} HS-SVD
&
Orthogonal expansion combined with SVD
&
{\small Very flat kernels and spectral accuracy are required}
&
$O(N^3)$
&
Excellent
&
Moderate
&
{\small Mathematically sophisticated; basis depends on kernel and geometry. }\\

\rowcolor{blue!20} WSVD
&
{\small Weighted SVD producing an orthonormal basis}
&
{\small General-purpose stabilization for arbitrary kernels}
&
$O(N^3)$
&
Very good
&
Moderate
&
{\small Dense SVD becomes expensive for large datasets.} \\

\rowcolor{blue!20} Greedy
&
{\small Adaptive selection of interpolation centers}
&
{\small Large datasets requiring reduced approximation spaces}
&
Typically $O(mN^2)$,
$m\ll N$
&
Indirectly good
&
Good
&
{\small Accuracy depends on the selection criterion; does not directly remove ill-conditioning. }\\

\rowcolor{blue!20} PU
&
{\small Domain decomposition into local interpolation problems}
&
{\small Very large datasets and high-dimensional applications}
&
$\displaystyle\sum_j O(N_j^3)$
&
Good
&
Excellent
&
{\small Requires tuning of patches, overlap and weight functions. }\\

\bottomrule
\end{tabular}
}
\end{table*}

The {\em $f$-greedy} \cite{DeMarchi2, Muller2} is a data-dependent method that selects the next point where the current interpolation residual {{reaches}} the maximum value, so this is highly effective for functions with localized features {{such as}} sharp gradients. 
The {\em $f/p$-greedy} \cite{Muller2} is a hybrid selection strategy designed to maximize the benefits of both data-independent ($p$-greedy) and data-dependent ($f$-greedy) methods. 
The {\em $\beta$-greedy} \cite{Wenzel} algorithm is a unified theoretical framework that encompasses the entire family of greedy selection strategies for RBF interpolation; it is formally introduced to study how different selection criteria affect the convergence rate and numerical stability of the interpolant.


We recall that the {\it Lebesgue constant}, $\Lambda_{N}$, is an indicator that measures the propagation of the interpolation error, indeed, we have
$$\|f-P_{f}\|_\infty \le \Lambda_N E^*_{N}(f),$$
where  $\|f-P_{f}\|_\infty =\max\{\mod{f(\vv x)-P_f(\vv x)}$, $\vv x\in D\}$, and  as usual $E^*_{N}(f)$ is the {{sup-norm}} error of the best polynomial approximation of $f$ of degree at most $N$.
It grows as the number of data points $N$ increases, but high-quality point sets $X$, such as those generated by the $p$-greedy algorithm, exhibit slow growth of the Lebesgue constant; see \cite{DeMarchi2} for details. We {also} note that $p$-greedy or their generalization known as $\beta$-greedy strategies ensure a good separation distance among points in $X$, preventing in this way a fast {{growth}} of the condition number of the system \en{RBFint}, see \cite{Wenzel, Santin} for details.

\begin{figure}[!h]
\centering\includegraphics[scale=0.5]{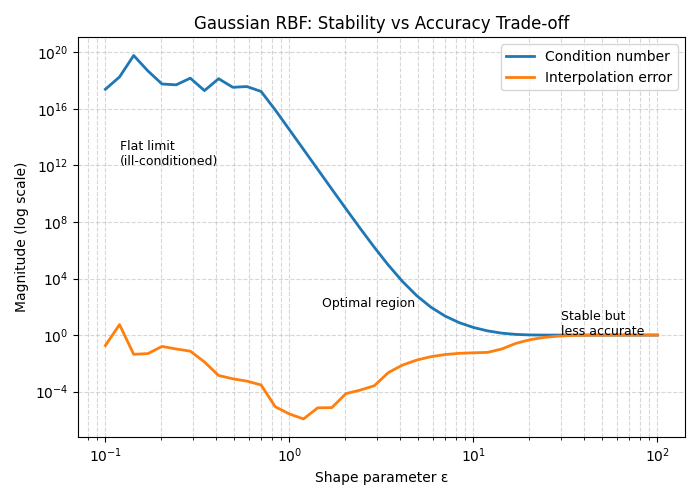}
\caption{The trade-off principle or the Shape Parameter Dilemma. {{Conditioning Vs Accuracy in the interpolation of the Franke function. In red the RMSE, in blue the condition number in 2-norm. In the \apap Optimal region\apch\ we have a good balance between accuracy and stability.  }} 
 }\label{fig:ShapeParamDilemma}
\end{figure}

A similar result for greedy strategies is obtained in \cite{Larsson}, where the statistical concepts of Design of Experiments are used to improve the numerical stability of the interpolation matrix and the interpolation error, {{thus}} mitigating the {\em Shape Parameter Dilemma}. It is worth noting that this yields a point distribution similar to Chebyshev or Leja points,
which are denser near the boundaries of the domain to avoid the {\em Runge effect}; moreover, a similar approach is used in the solution of partial differential equations. {{ 
The Shape Parameter Dilemma strongly affects both the effectiveness and the accuracy of an RBF interpolant; Figure~\ref{fig:ShapeParamDilemma} shows a trade-off in the shape parameter choice for Gaussian RBF. Both the flat limit region and the tiny variance region must be avoided due to ill-conditioning or scarce accuracy, respectively. The ``safe" choice is indicated in the figure as \apap Optimal region" (the central region) where the condition number is not exploding and the interpolation error is satisfactory.}}

{{We conclude the section with some comments on {\it Partition of Unity (PU)} techniques \cite{Babuska, Wendland}, as a \apap trait d'union" with the next section. Indeed, PU techniques are capable of stabilizing the problem \en{RBFint} but also providing a more efficient solution of RBF interpolation by transforming a single, often ill-conditioned global problem into a series of smaller, numerically manageable local subproblems.}}
The main components of these techniques are: {\bf i)} the reduction of ill-conditioning by dividing the domain into overlapping subdomains (patches), where small local matrices have to be inverted with a consequent reduction of numerical instability; {\bf ii)} local management of the problem by adapting the solution procedure to each patch, for example, by using local shape parameters; {\bf iii)} smooth blending by compactly supported weight functions to ensure a smooth, continuous transition between patches without introducing artificial oscillations. 
After its introduction, PU has been used in several concrete applications, such as, for example, 3D surface reconstruction from laser scans \cite{Ohtake}, Interpolation of Earth’s topography data \cite{Cavoretto}, and smooth transfer of aerodynamic loads from a dense fluid mesh to a sparser structural wing mesh \cite{Beckert}.

{In Tables \ref{tab:comparison_stabilization} is reported a synthesis of the behaviors of the considered stabilization methods. }

\section{Computational strategies: fast and compression}
\label{sec3}
We premise that any strategy to reduce the computational cost of the RBF interpolation should be implemented with special care, especially when these strategies are based on an approximation of the RBF matrix; in fact, any small modifications of the RBF matrix can produce large perturbations in the solution of \en{RBFint} as a consequence of its ill conditioning. {{In the following sections we will analyze in detail some of these techniques whose characteristics are summarized in the Table \ref{tab:comparison_acceleration}.}}

\subsection{Classical approaches}
The {\it Fast Multipole Method (FMM)} is exactly of this type; it is used to accelerate the evaluation of the RBF interpolant in the solution of the linear system \en{RBFint}. 
This is a dense linear system, so direct methods (like LU decomposition) require $O(N^3)$ operations, the Conjugate Gradient Method (CGM) converges in at most $N$ iterations, because the matrix is symmetric. A single matrix-vector multiplication (in iterative methods such as CGM or GMRES or BiCGSTAB) requires $O(N^2)$ operations. So, reducing the computational cost of matrix-vector multiplication is fundamental {{to}} obtaining efficient iterative methods {{to solve}} \en{RBFint}.
The FMM reduces the cost of matrix-vector multiplication in iterative methods by approximating the action of the RBF matrix. This approximation uses ideas from classical scattering theory \cite{Greengard} and can compute matrix-vector multiplication in $O(N\log N)$ operations. In \cite{DeMarchiE}, a fast algorithm {{for}} matrix-vector multiplication is obtained for the inverse multiquadric; the computational cost of this algorithm has an upper bound of $O(N\,\log N )$.
In \cite{Fong}, the Chebyshev interpolation is used to replace the analytical expansion of the kernel  $\phi$. 
In \cite{Ying}, the analytical expansion of the kernel $\phi$ is calculated in terms of equivalent densities on surfaces enclosing the domain $D$; this approach has general validity and allows the use of different kernels $\phi$. In \cite{Gumerov}, the acceleration of the iterative solver for RBFs using FMM has been provided.
{In \cite{Giacomini}, using a technique similar to FMM and a local low-rank representation of the interpolation matrix, a preconditioning matrix is obtained to stabilize the numerical solution of \en{RBFint}. In \cite{Egidi}, a preconditioning strategy for data interpolation through RBFs is devised by exploiting spectral techniques, such as the Jacobi method, for the reconstruction of the approximated eigenvalues of the original interpolation system.}

A similar method that allows both efficient matrix-vector multiplication and an explicit representation of the matrix, which can be profitably used to construct a preconditioner, is the \textit{hierarchical representation} \cite{Armin}; in this paper, {{different approaches are analyzed to obtain low-rank and block-wise approximations for various RBFs.}} Among the various approaches considered, a Taylor-based representation undoubtedly benefits from its generality, which {{guaranties}} its applicability to any conditionally positive definite RBF.

High-performance computing software has been developed to implement {{FMM,}} such as, for example, \cite{Yokota}  demonstrates {{the effectiveness}} of RBF interpolation for datasets with several million points using parallel FMM and the PETSc library; \cite{Owen} provides a computer procedure in the Rust language \cite{Rust}, combining FMM and domain decomposition for datasets with {{more than a}} million points. 

\subsection{Samplets: a new compression strategy}

\begin{figure}[!h]
\centering
\includegraphics[scale=0.6]{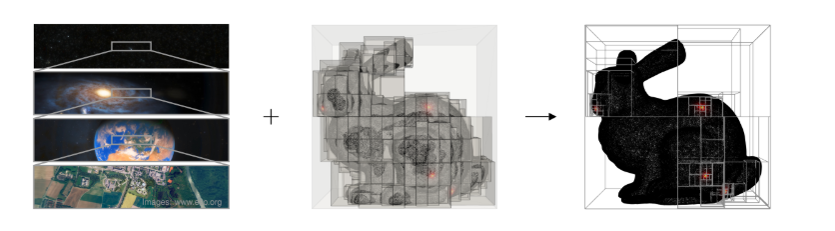}
\caption{Left: multiresolution of the Earth from the space. Center: the bunny overlayed by samplets. Right: the sum shows how the multiresolution of samplets works. Image courtesely provided by Sara Avesani - Universit\`a della Svizzera Italiana.} \label{fig:Samplets}
\end{figure}

A similar result can be obtained by {\em samplets} \cite{Harbrecht}, which are a specific class of multiscale basis functions designed to compress the dense matrices that arise in large-scale problems. 
In particular, they have been introduced as a specialized alternative to standard wavelets, specifically tailored for scattered data and integral operators. Samplets are a specialized multiresolution analysis (MRA) technique designed for data-adapted decomposition of localized signed measures, frequently used to handle scattered (unstructured) data in arbitrary dimensions. They represent a generalization of the Tausch-White multi-wavelet framework to discrete settings, enabling data compression, feature detection, and numerical approximation  \cite{Harbrecht}, and we may construct them using the usual vanishing moment framework. A visual representation on a multiresolution by samplets is given in Figure \ref{fig:Samplets}.

This abstract setting allows us to construct them on very general domains (hence data sets).
By \apap thresholding" {{(setting small values to zero),}} the dense matrix becomes sparse; in this way, {{matrix-vector multiplication requires $O(N\log N)$ operations that can be used profitably}} in the iterative solution of \en{RBFint}. The samplets framework has been extended to efficiently compute complex algebraic operations on the RBF matrix; see \cite{Harbrecht2} for details. In \cite{Baroli}, a study of sample properties for functions in $N_\phi(D)$ and some applications of RBF in sparse data interpolation and surface reconstruction are reported. 

\subsection{The Nystr\"{o}m method}
The {\em Nystr\"{o}m method} provides another approximation technique to reduce the computational cost of \en{RBFint}. In this technique, the kernel matrix $\Phi$ is approximated using only $m$ columns, where $m\ll N$. The selection of these landmarks from the dataset $X$ is usually {{performed}} by random sampling, clustering of K-means, or more advanced techniques such as the {\em Leverage Score sampling} \cite{Drineas}. The eigendecomposition of the selected $m$ columns provides a low-rank representation of the kernel matrix
\begin{equation}
\Phi\approx\tilde\Phi=UEU^t,
\label{nyst}
\end{equation}
where $U\in R^{N\times m}$, $E\in R^{m\times m}$. From \en{nyst}, the Woodbury formula can be used to solve the linear system \en{RBFint} efficiently.
Statistical analysis and error bounds for the method are provided in \cite{Drineas2}; moreover, a detailed analysis of the strategy to select better landmark points is proposed in \cite{Zhang}.

A similar result is obtained  by the {\em Structured Kernel Interpolation} (SKI) \cite{Wilson} that approximates the kernel matrix $\Phi$ using an interpolation on a regular grid, in this way 
\begin{equation}
\Phi\approx\tilde\Phi=WKW^t,
\label{ski}
\end{equation}
where {{$W\in R^{N\times M}$}} is a sparse matrix containing the interpolation weights and $K$ has a Toeplitz structure; {{therefore,}} the solution of \en{RBFint} can be obtained by iterative methods {{with a computational cost nearly linear.}} In \cite{Gardner}, it is proposed to improve the efficiency of SKI by a decomposition of the kernel into products of lower-dimensional structures. Moreover, \cite{Moreno} has provided a rigorous theoretical error bound for SKI, confirming that the approximation error decays cubically relative to the grid density when using cubic interpolation.

{In Table \ref{tab:comparison_acceleration} is reported a synthesis of the behaviors of the acceleration methods considered.}

\subsection{Other approaches}
The conceptual basis of {\em Random Fourier Feature} (RFF) \cite{Rahimi} is the integral Bochner Theorem, which is used to approximate the kernel $\phi$ by a finite number of randomized trigonometric features and the explicit mapping of the data to a low-dimensional Euclidean space. In this way, the solution of \en{RBFint} is obtained as a least squares problem in a reduced number of variables. The theoretical analysis for the uniform convergence of random features is provided in \cite{Rahimi2}. Several authors have studied this method in detail, analyzing uniform error bounds \cite{Sutherland}, {{the quality of the approximation,}} the statistical efficiency \cite{Sriperumbudur}, and the risk of learning with RFF for different loss functions \cite{Li}.

Among the various approximation techniques for {{the}} RBF matrix $\Phi$, it is worth mentioning also the so-called {\em Adaptive Cross Approximation} (ACA), initially proposed in the context of the solution of integral equations \cite{Bebendorf} and then generalized for RBF interpolation matrices \cite{Grzhibovskis}. This algorithm constructs low-rank approximations of a given matrix by selecting rows and columns, calculating the residual, and updating the matrices until a tolerance threshold is met.
ACA does not use analytical expansions of the kernel function $\phi$, so, in principle, it can work with different kernels even if it is developed for matrices generated by sampling smooth functions.  


\begin{table*}[ht]
\caption{Comparison of the main acceleration techniques for kernel-based interpolation and machine learning.}
\label{tab:comparison_acceleration}
{\footnotesize
\!\!\!\!\!\!\!\!\!\!\!\!\!\!\!\!\!\!\!\!\!\!\!\begin{tabular}{p{1.6cm}p{2.0cm}p{1.5cm}p{1.7cm}p{2.8cm}p{2.8cm}}
\toprule
\textbf{Method} &
\textbf{Cost} &
\textbf{Memory} &
\textbf{Approx.} &
\textbf{Use for} &
\textbf{Impl. issues} \\
\midrule

\rowcolor{blue!20} FMM
&
$O(N)$--$O(N\log N)$
&
$O(N)$
&
Very high
&
{\small Long-range kernel interactions, RBF interpolation, $N$-body problems}
&
{\small Kernel-dependent expansions, sophisticated implementation, difficult for arbitrary kernels. }\\

\rowcolor{blue!20} $\mathcal H$-matrices
&
$O(N\log N)$
&
$O(N\log N)$
&
High
&
{\small General kernel matrices arising from PDEs and interpolation}
&
{\small Requires hierarchical partitioning and low-rank block compression; many tuning parameters.} \\

\rowcolor{blue!20} Samplets
&
$\approx O(N)$
&
$O(N)$
&
High
&
{\small Multiscale kernel approximation on manifolds and scattered data}
&
{\small Construction of multiscale basis is nontrivial; mainly developed for specific kernels and geometries.} \\

\rowcolor{blue!20} Nyström
&
$\scriptstyle O(Nm+m^3)$, $m\ll N$
&
$O(Nm)$
&
{\small Good, depending on landmark selection}
&
{\small Large kernel matrices, Gaussian processes, kernel machines}
&
{\small Accuracy strongly depends on sampling strategy and number of landmarks. }\\

\rowcolor{blue!20} RFF
&
$O(ND)$
&
$O(ND)$
&
Moderate--High
&
{\small Shift-invariant kernels, large-scale learning}
&
{\small Requires kernel stationarity; accuracy depends on the number $D$ of random features.} \\

\rowcolor{blue!20} SKI
&
$\scriptstyle \approx O(N+M\log M)$
&
$O(N+M)$
&
High
&
{\small Gaussian processes and interpolation on structured grids}
&
{\small Requires interpolation onto inducing grids; most effective for tensor-product or nearly structured domains.} \\

\bottomrule
\end{tabular}
}
\end{table*}


\section{Conclusions and future perspectives}
\label{sec4}

Radial Basis Function approximation has evolved from a classical meshfree
interpolation technique into a mature computational framework with applications
ranging from scattered data approximation and the numerical solution of partial
differential equations to machine learning, image registration, mesh morphing,
and scientific computing \cite{Kansa,Biancolini,Boer,Rohr,Broomhead,Moody,Poggio,Boser,
Fornberg3,Li2,Schaback3,Haykin,Cristianini,Bandiziol24}. The stabilization
techniques discussed in this review, including RBF--QR
\cite{Fornberg11,Kormann19,Fornberg1,Fornberg2,Fasshauer2}, HS--SVD
\cite{Cavoretto15bb,FasshauerMcCourt}, WSVD bases \cite{DeMarchi13}, greedy
algorithms \cite{
DeMarchi,DeMarchi2,Wenzel,Muller2,Schaback2,Santin,Larsson},
and Partition of Unity methods \cite{Babuska,Wendland,Ohtake,Cavoretto,Beckert},
have considerably reduced the numerical difficulties associated with
ill-conditioning and the shape parameter dilemma. Similarly, fast computational
techniques such as the Fast Multipole Method
\cite{Greengard,DeMarchiE,Fong,Ying,Gumerov,Giacomini}, hierarchical matrices
\cite{Armin}, samplets \cite{Harbrecht,Harbrecht2,Baroli}, Nystr\"om
approximation \cite{Drineas,Drineas2,Zhang,Williams}, Random Fourier Features
\cite{Rahimi,Rahimi2,Sutherland,Sriperumbudur,Li}, and Structured Kernel
Interpolation \cite{Wilson,Gardner,Moreno} have substantially improved the
scalability of kernel-based algorithms, making RBF methods increasingly
attractive for large-scale scientific computing.

Despite these remarkable developments, several important challenges remain
open. One of the most significant is the automatic selection of the shape
parameter. Although numerous heuristic and optimization-based strategies have
been proposed, a general theory capable of balancing numerical stability,
approximation accuracy, and computational efficiency is still missing. In
particular, adaptive strategies that combine local error indicators, point
distributions, and optimization techniques appear to represent one of the most
promising research directions. Such approaches could naturally complement
greedy center selection \cite{Wenzel,Schaback2,Santin,Larsson} and local
Partition of Unity constructions \cite{Babuska,Wendland,Cavoretto}.

Another important challenge concerns the development of RBF algorithms for
next-generation high-performance computing architectures. While parallel
implementations based on Fast Multipole Methods, hierarchical representations,
and domain decomposition have already demonstrated excellent performance
\cite{Gumerov,Armin,Yokota,Owen}, the transition towards exascale computing
requires algorithms capable of exploiting heterogeneous architectures composed
of multicore CPUs, GPUs, high-bandwidth memories, and distributed-memory
systems. In this context, communication-avoiding algorithms, asynchronous
implementations, scalable preconditioners, and fault-tolerant solvers will
probably become key ingredients for future developments.

Closely related to this aspect is the increasing importance of GPU-oriented
algorithms. The dense linear algebra underlying global RBF interpolation can
benefit greatly from massive parallelism, but efficient implementations require
a careful redesign of kernel evaluations, hierarchical matrix operations,
adaptive tree structures, and iterative solvers in order to maximize memory
locality and computational throughput. Existing parallel RBF software and
special-purpose hardware \cite{Yokota,Hines,Crean,Bird,CM1K,SKhynix} provide
useful starting points, but the combination of stabilized bases with GPU and
heterogeneous implementations remains largely unexplored.

A further rapidly growing area concerns uncertainty quantification. Many
applications involve noisy measurements, incomplete datasets, or uncertain
model parameters. Since kernel approximation is closely related to Gaussian
process regression and probabilistic learning, integrating RBF approximation
with Bayesian methodologies, confidence estimation, stochastic inverse
problems, and uncertainty propagation represents a natural extension of the
deterministic framework presented in this survey. An important theoretical
question is how numerical stabilization, kernel compression, and uncertainty
estimates interact when the kernel matrix is ill-conditioned or only
approximately represented.

Equally important is the recent interest in kernel learning. Classical RBF
methods usually assume that both the kernel family and its parameters are fixed
a priori. Modern machine learning instead suggests learning kernel parameters
directly from the data, combining multiple kernels, or constructing
problem-adapted kernel representations. Such approaches have the potential to
combine the mathematical rigor and interpretability of classical kernel
approximation with the flexibility of data-driven statistical learning. 
In
this setting, scalable approximations such as Nystr\"om methods, Random Fourier
Features, and Structured Kernel Interpolation
\cite{Drineas2,Wilson,Gardner,Rahimi,Li,Williams} may play a particularly
important role.

Finally, the interaction between RBF approximation and modern scientific
machine learning is expected to become increasingly important. Besides their
well-established role in Radial Basis Function Neural Networks
\cite{Broomhead,Moody,Poggio,Haykin} and Support Vector Machines
\cite{Boser,Cristianini}, kernel methods are beginning to interact with
physics-informed neural networks, operator learning, and neural operators
\cite{86,87,89,90}. Hybrid methodologies combining stable kernel approximation
with deep neural architectures may provide accurate surrogate models,
interpretable learning algorithms, and reliable uncertainty estimates for
complex scientific simulations.

Overall, the future of RBF approximation will likely be driven by the
integration of three complementary research directions: improved numerical
stability, scalable high-performance implementations, and adaptive data-driven
methodologies. The combination of stabilization techniques, compression
strategies, heterogeneous computing architectures, uncertainty quantification,
and modern machine learning appears to be one of the most promising avenues for
the next generation of kernel-based numerical methods. Consequently, RBF
approximation continues to represent a vibrant research area at the
intersection of approximation theory, numerical analysis, scientific computing,
and artificial intelligence.

\vskip 0.1in
{\bf Acknowledgments.} {{We thank the anonymous reviewers who allowed us to improve the manuscript with their fruitful comments and recommendations.}} We acknowledge the INdAM-GNCS group and the topical group on ``Approximation Theory and Applications'' of the Italian Mathematical Union and the topical group ANA\&A of SIMAI.


\begin{thebibliography}{99}
\bibliographystyle{abbrv}
\bibitem{hardy} R.L. Hardy, Multiquadric equations of topography and other irregular surfaces. \emph{Journal of Geophysical Research,} {\bf 76} (1971) 1905-1915.
\bibitem{Meinguet} J. Meinguet,  Multivariate interpolation at arbitrary points made simple. \emph{Zeitschrift für angewandte Mathematik und Physik (ZAMP),} {\bf 30} (1979) 292-304.
\bibitem{Franke} R. Franke, Scattered data interpolation: tests of some methods. \emph{Mathematics of Computation,} {\bf 38} (1982), 181-200.
\bibitem{Micchelli} C.A. Micchelli,   Interpolation of scattered data: Distance matrices and conditionally positive definite functions. \emph{Constructive Approximation,} {\bf 2} (1986), 11-22.
\bibitem{Duchon} J. Duchon. Splines minimizing rotation-invariant semi-norms in Sobolev spaces. In: \emph{W. Schempp, K. Zeller (eds) Constructive Theory of Functions of Several Variables. Lecture Notes in Mathematics,} vol 571, 1977.

\bibitem{Kansa} E.J. Kansa, Multiquadrics - A scattered data approximation scheme with applications to computational fluid-dynamics-II. \emph{Computers \& Mathematics with Applications,} {\bf 19} (1990) 147-161.
\bibitem{Biancolini} M.E. Biancolini,  Fast Radial Basis Functions for Engineering Applications,  2018. Springer International Publishing.
\bibitem{Boer}
A. de Boer, M.S. van der Schoot, H. Bijl, Mesh deformation based on radial basis function interpolation. \emph{Computers \& Structures,} {\bf 85} (2007) 784-795.
\bibitem{Rohr}
K. Rohr, H.S. Stiehl, R. Sprengel, T.M. Buzug, J. Weese, M.H. Kuhn, Landmark-based elastic registration using approximating thin-plate splines. \emph{IEEE Transactions on Medical Imaging,} {\bf 20} (2001) 526-534.
\bibitem{Broomhead} D.S. Broomhead, D. Lowe,  Multivariable functional interpolation and adaptive networks. \emph{Complex Systems,} {\bf 2} (1988) 321-355.
\bibitem{Moody} J. Moody, C.J. Darken, Fast learning in networks of locally-tuned processing units. \emph{Neural Computation,} {\bf 1} (1989) 281-294.

\bibitem{Poggio} T. Poggio, F. Girosi, Networks for approximation and learning. \emph{Proceedings of the IEEE,} {\bf 78} (1990) 1481-1497.

\bibitem{Boser}
B.E. Boser, I.M. Guyon, V.N. Vapnik, A training algorithm for optimal margin classifiers. \emph{COLT '92: Proceedings of the fifth annual workshop on Computational learning theory,}
(1992) 144-152.

\bibitem{Fasshauer} G.E. Fasshauer, Meshfree Approximation Methods with MATLAB. \emph{Interdisciplinary Mathematical Sciences,} Vol. 6, 2007, World Scientific Publishers.
\bibitem{Fornberg} B. Fornberg, G. Wright, Stable computation of multiquadric interpolants for all values of the shape parameter. 
\emph{Computers \& Mathematics with Applications}, {\bf 48} (2004) 853-867


\bibitem{DeMarchi} S. De Marchi, R. Schaback, H. Wendland, Near-optimal search of centers for RBF prediction models. \emph{Computers \& Mathematics with Applications,} {\bf 50} (2005) 369–380.



\bibitem{DeMarchi2} S. De Marchi, Geometric greedy and greedy points for RBF interpolation. \emph{Proceedings of the International Conference on Computational and Mathematical Methods in Science and Engineering,} 30 June, 1–3 July 2009.

\bibitem{Wenzel} T.  Wenzel, G. Santin, B. Haasdonk, Analysis of Target Data-Dependent Greedy Kernel Algorithms. \emph{Constructive Approximation,} {\bf 55} (2022) 871-911.

\bibitem{Muller2} S. Müller, Komplexität und Stabilität von kernbasierten Rekonstruktionsmethoden. PhD thesis, University of Göttingen, (2009)

\bibitem{Eijkhout} V. Eijkhout, E. Chow, R. van de Geijn, Introduction to High Performance Scientific Computing, CC BY 3.0, 3rd edition 2020
\bibitem{Narcowich1} F.J. Narcowich, J.D. Ward, Norm of inverses and condition numbers for matrices associated with scattered data. \emph{J. Approx. Theory,} {\bf 64} (1991) 69–94.

\bibitem{Narcowich2} F.J. Narcowich, J.D. Ward, Norms of inverses for matrices associated with scattered data. \emph{In: Curves and Surfaces, eds. P.J. Laurent, A. Le Méhauté and L.L. Schumaker,} (Academic Press, Boston, 1991) pp. 341–348.

\bibitem{Narcowich3} F.J. Narcowich, J.D. Ward, Norm estimates for the inverses of a general class of scattered-data radial-function interpolation matrices. \emph{J. Approx. Theory,} {\bf 69} (1992) 84–109.

\bibitem{Ball1} K. Ball, N. Sivakumar, J.D. Ward, On the sensitivity of radial basis interpolation to minimal data separation distance. \emph{Constr. Approx.,} {\bf 8} (1992) 401–426.

\bibitem{Ball2} K. Ball, Eigenvalues of Euclidean distance matrices. \emph{J. Approx. Theory,} {\bf 68} (1992) 74-82.

\bibitem{Braun} M.L. Braun, Accurate Error Bounds for the Eigenvalues of the Kernel Matrix. \emph{Journal of Machine Learning Research,} {\bf 7} (2006) 2303-2328.

\bibitem{Schaback} R. Schaback,  Error estimates and condition numbers for radial basis function interpolation. \emph{Adv. Comput. Math.,} {\bf 3} (1995) 251–264.
\bibitem{DeMar} S. De Marchi, R. Schaback, H. Wendland, Near-Optimal Data-independent Point Locations for Radial Basis Function Interpolation. \emph{Adv. Comput. Math.,} {\bf 23} (2005) 317-330.

\bibitem{FoPi07}{B.~Fornberg and C.~Piret},
A stable algorithm for flat radial basis functions on a sphere. \emph{SIAM J. Sci. Comput.,} \textbf{30} (2007), pp.~{60--80}.

\bibitem{Fornberg11} {B.~Fornberg, E.~Larsson, and N.~Flyer},
Stable computations with {G}aussian radial basis functions. \emph{SIAM J. Sci. Comput.,} \textbf{33} (2011), pp.~869--892.

\bibitem{FasshauerMcCourt12} G.E. Fasshauer and M. Mc Court,  Stable evaluation of Gaussian radial basis
function interpolants. \emph{SIAM J. Sci. Comput.,} {\bf 34} (2012), pp. A737–A762.

\bibitem{Kormann19} {K.~Kormann, C.~Lasser, and A.~Yurova},
Stable interpolation with isotropic and anisotropic Gaussians using Hermite generating function. \emph{SIAM J. Sci. Comput.,} \textbf{41} (2019), pp.~A3839--A3859.ù

\bibitem{Cavoretto15bb}
{R. Cavoretto and G.E. Fasshauer and M. McCourt}, {An introduction to the Hilbert-Schmidt SVD using iterated Brownian bridge kernels}.
\emph{Numer. Algorithms,} \textbf{68} (2015), pp. {393--422}

\bibitem{FasshauerMcCourt} G.E. Fasshauer and M. Mc Court, Kernel-based Approximation Methods with MATLAB. \emph{Interdisciplinary Mathematical Sciences,} Vol. 19, 2015, World Scientific Publishers.
\bibitem{DeMarchi13}
{S. {De Marchi} and G. Santin},
A new stable basis for radial basis function interpolation. \emph{J. Comput. Appl. Math.,} \textbf{253} (2013),pp. {1-13}.
\bibitem{Fornberg1} B. Fornberg and C. Piret, A stable algorithm for flat radial basis functions on a sphere, \emph{SIAM Journal on Scientific Computing,} {\bf 30} (2007) 60–80.

\bibitem{Fornberg2} B. Fornberg, E. Larsson, N. Flyer, Stable computations with gaussian radial basis functions. \emph{SIAM Journal on Scientific Computing,} {\bf 33} (2011) 869–892.

\bibitem{Fasshauer2} G.E. Fasshauer, M.J. McCourt, Stable evaluation of Gaussian radial basis function interpolants. \emph{SIAM Journal on Scientific Computing,} {\bf 34} (2012) A737–A762.
\bibitem{Schaback2} R. Schaback, H. Wendland, Adaptive greedy techniques for approximate solution of large RBF systems. \emph{Numerical Algorithms,} {\bf 24} (2000) 239–254 .
\bibitem{Santin} G. Santin, B. Haasdonk, Convergence rate of the data-independent P-greedy algorithm in kernel-based approximation. \emph{Dolomites Research Notes on Approximation,} {\bf 10} (2017) 68-78.
\bibitem{Larsson} E. Larsson, R. Schaback, A Design of Experiments Approach to Point Selection in Radial Basis Function Methods, \emph{Engineering Analysis with Boundary Elements,}  {\bf 31} (2007) 64–72.
\bibitem{Babuska} I. Babu\v{s}ka, J.M. Melenk, The partition of unity method. \emph{International Journal for Numerical Methods in Engineering,} {\bf 40} (1997) 727-758.

\bibitem{Wendland} H. Wendland, Fast evaluation of radial basis functions: Methods based on partitions of unity. \emph{Approximation Theory X: Wavelets, Splines, and Applications,} 2002.

\bibitem{Ohtake} Y. Ohtake, A. Belyaev, M. Alexa, G. Turk, H.-P. Seidel, Multi-level partition of unity implicits. \emph{ACM Transactions on Graphics,} {\bf 22} (2003) 463-470.

\bibitem{Cavoretto} R. Cavoretto, A. De Rossi, A meshless interpolation algorithm using a cell-based searching procedure. \emph{Computers \& Mathematics with Applications,} {\bf 67} (2014) 1024-1038.

\bibitem{Beckert} A. Beckert, H. Wendland, Multivariate interpolation for fluid-structure-interaction using radial basis functions. \emph{Aerospace Science and Technology,} {\bf 5} (2001) 125-134.

\bibitem{Greengard} L.F. Greengard, V. I. Rokhlin, A Fast Algorithm for Particle Simulation. \emph{Journal of Computational Physics,} {\bf 73} (2001) 325-348.

\bibitem{DeMarchiE} S. De Marchi, N. Egidi, J. Giacomini, P.
Maponi, A. Perticarini,  Computational issues by interpolating
with inverse multiquadrics: a solution. \emph{Dolomites Research Notes
on Approximation,} {\bf 15} (2022) 56-64.
\bibitem{Fong} W. Fong, E. Darve, The black-box fast multipole method. \emph{Journal of Computational Physics,} {\bf 228} (2009) 8712–8725.

\bibitem{Ying} L. Ying, G. Biros, D. Zorin, A kernel-independent adaptive fast multipole algorithm
in two and three dimensions. \emph{Journal of Computational Physics,} {\bf 196} (2004) 591–626.

\bibitem{Gumerov} N. Gumerov, R. Duraiswami, Fast Radial Basis Function Interpolation via Preconditioned Krylov Iteration. \emph{SIAM Journal on Scientific Computing,} {\bf 29} (2007) 1876-1899.





\bibitem{Giacomini} J. Giacomini,  RBFs preconditioning via Fourier decomposition method. \emph{AIP Conference Proceedings} {\bf 3094}(1) (2024), 320005.

\bibitem{Egidi} N. Egidi, J. Giacomini, P.
Maponi,  Preconditioning strategies for RBF interpolation. \emph{Yaroslav D. Sergeyev, Dmitri E. Kvasov, and Annabella Astorino, editors, Numerical Computations: Theory and Algorithms,} Springer Nature, Switzerland, (2025) 246–253.

\bibitem{Armin} A. Iske, S. Le Borne, and M. Wende, Hierarchical matrix
approximation for kernel-based scattered data interpolation. \emph{SIAM
journal on scientific computing,} {\bf 39}(5) (2017) A2287–A2316.


\bibitem{Yokota} R. Yokota, L.A. Barba, M.G. Knepley, PetRBF - A parallel $O(N)$ algorithm for radial basis function interpolation with Gaussians. \emph{Computer Methods in Applied Mechanics and Engineering,} {\bf 199} (2010) 1793-1804.

\bibitem{Owen} D. Owen, ferreus\_rbf, Maptek Pty Ltd, 2025. Accessed March 12, 2026. 
\url{https://docs.rs/ferreus_rbf/latest/ferreus_rbf/#references}

\bibitem{Rust} The Rust Project Developers. \apap The Rust Programming Language." Accessed March 12, 2026. \url{https://www.rust-lang.org/}

\bibitem{Harbrecht} H. Harbrecht, M. Multerer, Samplets: Construction and scattered data compression. \emph{Journal of Computational Physics,} {\bf 471} (2022) 111616.

\bibitem{Harbrecht2} H. Harbrecht, M. Multerer, O. Schenk, Ch. Schwab, Multiresolution kernel matrix algebra. \emph{Numerische Mathematik,} {\bf 156} (2024) 1085–1114.

\bibitem{Baroli} D. Baroli, H. Harbrecht, M. Multerer, Samplet Basis Pursuit: Multiresolution Scattered Data Approximation With Sparsity Constraints. \emph{IEEE Transactions on Signal Processing,} {\bf 72} (2024) 1813-1823.
\bibitem{Drineas} P. Drineas, M.W. Mahoney, S. Muthukrishnan, Sampling algorithms for $\ell_2$ regression and applications. \emph{In Proceedings of the 17th Annual ACM-SIAM Symposium on Discrete Algorithms,}
pages 1127–1136, 2006.

\bibitem{Drineas2} P. Drineas, M.W. Mahoney, On the Nystr\"{o}m method for approximating a Gram matrix for improved kernel-based learning. \emph{Journal of Machine Learning Research,} {\bf 6} (2005) 2153–2175.

\bibitem{Zhang} K. Zhang, I.W. Tsang, J.T. Kwok, Improved Nyström low-rank approximation and error analysis. \emph{In: Proceedings of the 25th International Conference on Machine Learning,} Helsinki, Finland, 2008.

\bibitem{Wilson} A.G. Wilson, H. Nickisch, Hannes, Kernel interpolation for scalable structured Gaussian processes (KISS-GP). \emph{In: Proceedings of the 32nd International Conference on Machine Learning} - Volume 37, 2015, 1775–1784.
\bibitem{Gardner} J. Gardner, G. Pleiss, R. Wu, K. Weinberger, A. Wilson,  Product Kernel Interpolation for Scalable Gaussian Processes. \emph{In: Proceedings of the Twenty-First International Conference on Artificial Intelligence and Statistics,} 2018, 1407-1416.

\bibitem{Moreno} A. Moreno, J. Xiao, J. Mei, The Price of Linear Time: Error Analysis of Structured Kernel Interpolation. \emph{In: Proceedings of the Forty-Second International Conference on Machine Learning,} 2025.
\bibitem{Rahimi} A. Rahimi, B. Recht, Random Features for Large-Scale Kernel Machines, \emph{Advances in Neural Information Processing Systems,} {\bf 20} (2007) 1177--1184.

\bibitem{Rahimi2} A. Rahimi, B. Recht, Uniform approximation of functions with random bases. \emph{In: Proceedings of the46th Annual Allerton Conference on Communication, Control, and Computing,} 2008, 555--561.

\bibitem{Sutherland} D.J. Sutherland, J. Schneider, On the Error of Random Fourier Features. \emph{Proceedings of the 31st Conference on Uncertainty in Artificial Intelligence (UAI),} 2015, 862--871.
 
\bibitem{Sriperumbudur} B.K. Sriperumbudur, Z. Szab\'{o},  Optimal rates for Random Fourier features. \emph{Proceedings of the 29th International Conference on Neural Information Processing Systems} - Volume 1, 2015, 1144–1152.

\bibitem{Li} Z. Li, J. Ton, D. Oglic, D. Sejdinovic, Towards a Unified Analysis of Random Fourier Features. \emph{Journal of Machine Learning Research,} {\bf 22} (2021) 1--51.





\bibitem{Bebendorf} M. Bebendorf, Approximation of boundary element matrices. \emph{Numerische Mathematik,}  {\bf 86} (2000) 565–589.

\bibitem{Grzhibovskis} R. Grzhibovskis, M. Bambach, S. Rjasanow, G. Hirt, Adaptive cross-approximation for surface reconstruction using radial basis functions. \emph{Journal of Engineering Mathematics,} {\bf 62} (2008) 149–160.

\bibitem{Fornberg3} B. Fornberg, N. Flyer, A Primer on Radial Basis Functions with Applications to the Geosciences, 2015. Society for Industrial and Applied Mathematics.

\bibitem{Li2} J.-P. Li, Q.-H. Qin, Radial Basis Function Methods For Large-Scale Wave Propagation, 2021. Bentham Science publishers.

\bibitem{Schaback3} R. Schaback, H. Wendland, Kernel techniques: from machine learning to meshless methods. \emph{Acta numerica,} {\bf 15} (2006) 543-639

\bibitem{Haykin} S.S. Haykin, Neural Networks: A Comprehensive Foundation (2nd ed.), 1998. Prentice Hall.

\bibitem{Cristianini} N. Cristianini, J. Shawe-Taylor, An introduction to Support Vector Machines and other kernel-based learning methods, 2000. Cambridge University Press.

\bibitem{Bandiziol24} C. Bandiziol, S. De Marchi, Persistence symmetric kernels for classification: A comparative study.
\emph{Symmetry,} {\bf 16}(9) (2024), 1236.
\bibitem{Williams} C.K.I. Williams, M. Seeger, Using the Nystr\"{o}m method to speed up kernel machines. \emph{In: Advances in Neural Information Processing Systems,} (2001) 682-688.
\bibitem{Hines} T. Hines, RBF: A Python package for radial basis functions, 2022. GitHub.

\bibitem{Crean} J. Crean,  RadialBasisFunctions.jl: Radial basis functions for JuliaJ, 2024. GitHub.

\bibitem{Bird} T.-A. Bird, RadialBasisFiniteDifferences.jl: High-performance RBF-generated finite differences in Julia, 2024. GitHub.

\bibitem{CM1K} General Vision Inc. (2017). CM1K - NeuroMem Chip with 1024 Neurons Datasheet.

\bibitem{SKhynix} SK hynix Inc. (2022). SK hynix Develops PIM, Next-Generation AI Accelerator. SK hynix Newsroom

\bibitem{86}
M. Raissi, P. Perdikaris, G.E. Karniadakis,
Physics-informed neural networks: A deep learning framework for solving forward and inverse problems involving nonlinear partial differential equations.
\emph{Journal of Computational Physics},
{\bf 378} (2019), 686--707.

\bibitem{87}
L. Lu, P. Jin, G. Pang, Z. Zhang, G.E. Karniadakis,
Learning nonlinear operators via DeepONet based on the universal approximation theorem of operators. \emph{Nature Machine Intelligence},
{\bf 3} (2021), 218--229.

\bibitem{89}
Z. Li, N. Kovachki, K. Azizzadenesheli, B. Liu,
K. Bhattacharya, A. Stuart, A. Anandkumar,
Fourier Neural Operator for Parametric Partial Differential Equations.
In: \emph{International Conference on Learning Representations (ICLR)},
2021.

\bibitem{90}
N. Kovachki, S. Lanthaler, S. Mishra,
Neural Operator: Learning Maps Between Function Spaces. \emph{Acta Numerica},
{\bf 32} (2023), 1--97.

















\end{thebibliography}
\end{document}